\documentclass[12pt]{amsart}
\usepackage{mathrsfs}
\usepackage{amsmath,amsthm,amsfonts,amssymb,verbatim,eucal}

\numberwithin{equation}{section}

\newtheorem{thm}[equation]{Theorem}

\newtheorem{example}[equation]{Example}
\newtheorem{remark}[equation]{Remark}

\newenvironment{ex}{\begin{example}\rm}{\end{example}}

\newtheoremstyle{named}{}{}{\itshape}{}{\bfseries}{\!\!\!\!
 .}{.5em}{\thmnote{#3 }#1}
\theoremstyle{named}
\newtheorem*{namedthm}{}
\newcommand\quotient[2]{
        \mathchoice
            {% \displaystyle
                \, \text{\raise.45ex\hbox{$#1$}\big/\lower.45ex\hbox{$#2$}}\,%
            }
            {% \textstyle
                #1\,/\,#2
            }
            {% \scriptstyle
                #1\,/\,#2
            }
            {% \scriptscriptstyle  
                #1\,/\,#2
            }
    }

\DeclareMathOperator{\Ext}{Ext}
\DeclareMathOperator{\gr}{gr} 
\DeclareMathOperator{\Ima}{Im}

\DeclareMathOperator{\op}{op}

\newcommand{\NN}{\mathbb N}

\newcommand{\DOT}{\setlength{\unitlength}{1pt}\begin{picture}(2.5,2)
               (1,1)\put(2,3.5){\circle*{3}}\end{picture}}

\newcommand{\Z}{{\mathbb Z}}
\newcommand{\N}{{\mathbb N}}
\newcommand{\id}{\mbox{\rm id\,}}      
\newcommand{\Hom}{\mbox{\rm Hom\,}}

\renewcommand{\ker}{\mbox{\rm Ker\,}}

\newcommand{\ot}{\otimes}

\newcommand{\cH}{\mathcal{H}}
\newcommand{\CC}{\mathbb{C}}

\newcommand{\g}{{\mathfrak{g}}}
\newcommand{\del}{\partial}

\newcommand{\HH}{{\rm HH}}

\newcommand{\PBW}{Poincar\'e-Birkhoff-Witt}
\newcommand{\LH}{\text{LH}}
\newcommand{\GB}{\mathscr{G}}
\newcommand{\LM}{L\! M}                                   
\begin{document}
\begin{abstract}
{\sl }

We sample some Poincar\'e-Birkhoff-Witt theorems appearing in
mathematics. Along the way, we compare modern techniques used
to establish such results, for example, the Composition-Diamond 
Lemma, Gr\"obner basis theory, and the homological approaches of
Braverman and Gaitsgory and of Polishchuk and Positselski. We
discuss several contexts for 
PBW theorems and their applications, such as Drinfeld-Jimbo 
quantum groups, graded Hecke algebras, 
and symplectic reflection  and related algebras.

\end{abstract}
\title[Poincar\'e-Birkhoff-Witt theorems]
{Poincar\'e-Birkhoff-Witt theorems}
\date{April 25, 2014}
\author{Anne V.\ Shepler}
\address{Department of Mathematics, University of North Texas,
Denton, Texas 76203, USA}
\email{ashepler@unt.edu}
\author{Sarah  Witherspoon}
\address{Department of Mathematics\\Texas A\&M University\\
College Station, Texas 77843, USA}\email{sjw@math.tamu.edu}
\thanks{
This material is based upon work supported by the National
Science Foundation under Grant No.\ 0932078000, while the
second author was in residence at the Mathematical
Sciences Research Institute (MSRI) in Berkeley, California,
during the spring semester of 2013. 
The first author was partially supported by NSF grant
\#DMS-1101177.
The second author was partially supported by
NSF grant 
\#DMS-1101399.
}

\maketitle

%%%%%%%%%%%%%%%%%%%%%%%%%%%%%%%%%%%%%%%%%%%%%%%%%%%%%%%%%%%%%%
\section{Introduction}

In 1900, Poincar\'e~\cite{Poincare} 
published a fundamental result on Lie
algebras that would prove a powerful tool in representation theory:
A Lie algebra embeds into an associative algebra that behaves in
many ways like a polynomial ring. 
In 1937, Birkhoff~\cite{Birkhoff} 
and Witt~\cite{Witt} 
independently formulated and proved versions 
of the theorem that we use today, although neither author
cited Poincar\'e's earlier work.
The result was called the Birkhoff-Witt Theorem 
for years and then later the Poincar\'e-Witt Theorem 
(see Cartan and Eilenberg~\cite{CartanEilenberg}) before
Bourbaki~\cite{Bourbaki} prompted use of
its current name, the {\em Poincar\'e-Birkhoff-Witt Theorem}.

The original theorem on Lie algebras was greatly expanded over time by a
number of authors
to describe various algebras, 
especially those defined by quadratic-type
relations (including
Koszul rings over semisimple algebras).
\PBW\ theorems are often used as a springboard for investigating
the representation theory of algebras.  These theorems are used to
\begin{itemize}
\item reveal an algebra as a deformation of another, well-behaved
algebra,
\item posit a convenient basis (of ``monomials'') for an algebra, and
\item endow an algebra with  a canonical homogeneous (or graded) 
version.
\end{itemize}

In this survey, we sample some of the various 
Poincar\'e-Birkhoff-Witt theorems, applications,
and techniques used to date
for proving these results.
Our survey is not intended to be all-inclusive;
we instead seek to highlight a few of the more recent contributions and provide
a helpful resource for users of Poincar\'e-Birkhoff-Witt theorems, which
we henceforth refer to as {\em PBW theorems}. 

We begin with a quick review in Section~\ref{classical} 
of the original PBW 
Theorem for enveloping algebras of Lie algebras.
We next discuss PBW properties for quadratic algebras
in Section~\ref{homogeneous}, and
for Koszul algebras 
in particular,
before turning to arbitrary finitely generated
algebras in Section~\ref{sec:nonhomdef}.
We recall needed facts on
Hochschild cohomology
and algebraic deformation theory in Section~\ref{defHH},
and more 
background on Koszul algebras is given in Section~\ref{Koszul}.
Sections~\ref{BG}--\ref{diamond} outline techniques
for proving PBW results recently used in more general settings, some
by way of homological methods and others via the Composition-Diamond Lemma
(and Gr\"obner basis theory).
One inevitably is led to 
similar computations when applying any of
these techniques to specific algebras, but with different points of 
view. 
Homological approaches can help to organize computations and
may contain additional information, while approaches using
Gr\"obner basis theory are particularly well-suited for
computer computation.
We focus on some classes of algebras
in Sections~\ref{DJQG} and~\ref{SRA} of recent interest: 
Drinfeld-Jimbo quantum groups, Nichols algebras of diagonal type,
symplectic reflection algebras,
rational Cherednik algebras, and graded (Drinfeld) Hecke  algebras.
In Section~\ref{positivechar}, we mention 
applications in positive characteristic (including
algebras built on group actions in the modular
case) and other generalizations that mathematicians have only just begun to explore. 

We take all tensor products over an underlying field $k$ unless
otherwise indicated and assume all algebras are associative
$k$-algebras with unity.  Note that although we limit
discussions to finitely generated algebras over $k$ for simplicity,
many remarks extend to more general settings.

%%%%%%%%%%%%%%%%%%%%%%%%%%%%%%%%%%%%%%%%%%%%%%%%%%%%%%%%%%%%%%%%%%%%%
\section{Lie algebras and the classical PBW Theorem}\label{classical}
All PBW theorems harken back to
a classical theorem for universal enveloping algebras
of Lie algebras established independently
by Poincar\'e~\cite{Poincare}, Birkhoff~\cite{Birkhoff},
and Witt~\cite{Witt}.
In this section, we recall this original PBW 
theorem in order to set the stage for
other PBW theorems and properties; 
for comprehensive historical treatments, see \cite{Grivel,so-called}.

A finite dimensional {\em Lie algebra} is a finite dimensional
vector space $\g$ over a field $k$
together with a binary
operation $[ \ , \ ] : \g \times\g \rightarrow \g$ satisfying

(i)\ \ (antisymmetry) $\ \ \, [v,v]=0$ and

(ii)\, (Jacobi identity) $ \ [u,[v,w]]+[v,[w,u]]+[w,[u,v]] =0$

\noindent for all $u,v,w\in \g$.
Condition (i) implies $[v,w]=-[w,v]$ for all $v,w$ in $\mathfrak g$
(and is equivalent to this condition in all characteristics other than 2). 

The {\em universal enveloping algebra} $U(\g)$ of $\g$ is
the associative algebra generated by the vectors in $\g$ with relations
$vw-wv=[v,w]$ for all $v,w$ in $\g$, i.e.,
$$
  U(\g) = \quotient{T(\g)}{(v\ot w - w\ot v - [v,w] 
 : v,w\in \g), }
$$ 
where $T(\g)$ is the tensor algebra of the
vector space $\g$ over $k$. 
It can be defined by a universal property: $U(\g)$ is the
(unique up to isomorphism) associative algebra such that
any linear map $\phi$ from $\g$ to an associative algebra $A$
satisfying $[\phi(v),\phi(w)]=\phi([v,w])$ for all $v,w\in \g$
factors through $U(\g)$. 
(The bracket operation on an associative algebra $A$ is given 
by $[a,b]:= ab-ba$ for all $a,b\in A$.) 
As an algebra, $U({\mathfrak{g}})$ is filtered, under the assignment
of degree 1 to each vector in $\mathfrak{g}$. 

%%%%%%%%%%%%%%%%%%%
\begin{namedthm}[Original PBW Theorem]
A Lie algebra $\g$ embeds into its universal enveloping algebra
$U(\g)$, and the associated graded algebra of $U(\g)$ is isomorphic to 
$S(\g)$, the symmetric algebra on the vector space $\g$.
\end{namedthm}
%%%%%%%%%%%%%%%%%%%%%%%%%%%%%%%

Thus the original PBW Theorem 
compares a universal enveloping algebra $U(\g)$
to an algebra of (commutative) polynomials.
Since monomials form a $k$-basis for a polynomial algebra,
the original PBW theorem is often rephrased in terms of a 
{\em PBW basis} (with tensor signs between vectors dropped):

%%%%%%%%%%%%%%%%%%%
\begin{namedthm}[PBW Basis Theorem]
Let $v_1,\ldots,v_n$ be an ordered $k$-vector space basis of the Lie algebra $\g$.
Then
$
\{ v_{1}^{a_1}\cdots v_{n}^{a_n} : \  a_i \in \NN\}
$
is a $k$-basis of the universal enveloping algebra $U(\g)$.
\end{namedthm}
%%%%%%%%%%%%%%%%%%%%%%%%%%%%%5

\vspace{2ex}
\begin{ex}
The Lie algebra ${\mathfrak{sl}}_2(\CC)$ consists of $2\times 2$
matrices of trace 0 with entries in $\CC$
under the bracket operation on the 
associative algebra of all $2\times 2$ matrices.
The standard basis of ${\mathfrak{sl}}_2(\CC)$ is
\[
   e = \left(\begin{array}{cc} 0&1\\0&0\end{array}\right), \ \ 
   f = \left(\begin{array}{cc}0&0\\1&0\end{array}\right), \ \ 
   h = \left(\begin{array}{cc}1&0\\0&-1\end{array}\right),
\]
for which $[e,f]=h, \ [h,e]=2e,\ [h,f]=-2f$. Thus $U({\mathfrak{sl}}_2(\CC))$
is the associative $\CC$-algebra generated by three symbols
that we also denote by $e,f,h$ (abusing notation) subject to the relations
$ef-fe=h$, $\ he-eh=2e$, $\ hf-fh=-2f$.
It has $\CC$-basis 
$ \{ e^{a}h^{b}f^{c} : \, a,b,c\in\N\}$. 
\end{ex}
%%%%%%%%%%%%%%%%%%%%%%%%%%%%%%%%%%%%%%%%%%
\vspace{2ex}

Proofs of the original PBW Theorem vary (and by how much 
is open to interpretation).
The interested reader may wish to consult, for example, the texts~\cite{CartanEilenberg}, 
\cite{Dixmier}, 
\cite{Humphreys},
\cite{Jacobson}, 
and~\cite{Varadarajan}.
Jacobson~\cite{Jacobson41} proved a PBW theorem for restricted
enveloping algebras in positive characteristic.
Higgins~\cite{Higgins} gives references and a comprehensive PBW theorem over
more general ground rings. 
A PBW theorem for Lie superalgebras goes back to Milnor and Moore \cite{MilnorMoore}
(see also Kac~\cite{Kac}). 
Grivel's historical article~\cite{Grivel} includes further references on
generalizations to other ground rings, to Leibniz algebras, 
and to Weyl algebras. 
In Sections~\ref{BG} and \ref{diamond} below, we discuss two proof techniques particularly
well suited
to generalization:  a combinatorial approach through the 
Composition-Diamond Lemma 
and a homological approach through  algebraic deformation theory. 
First we lay some groundwork on quadratic algebras.

%%%%%%%%%%%%%%%%%%%%%%%%%%%%%%%%%%%%%%%%%%%%%%%%%%%%%%%%%%
%%%%%%%%%%%%%%%%%%%%%%%%%%%%%%%%%%%%%%%%%%%%%%%%%%%%%%%%%%
\section{Homogeneous quadratic algebras}\label{homogeneous}

Many authors have defined the notions of PBW algebra, 
PBW basis, PBW deformation, or PBW property 
in order to establish theorems like the
original PBW Theorem in more general settings. 
Let us compare a few of these concepts, beginning
in this section with those defined for {\em homogeneous} quadratic algebras.

\subsection*{Quadratic algebras}
Consider a finite dimensional vector space $V$ over $k$
with basis $v_1,\ldots, v_n$.
Let $T$ be its tensor algebra over $k$,  
i.e., the free $k$-algebra $k\langle v_1,\ldots, v_n\rangle $ generated by 
the $v_i$. Then $T$ is an $\NN$-graded $k$-algebra
with $$T^0=k,\ T^1=V,\ T^2=V\otimes V,\ T^3=V\otimes V\otimes V,
\text{ etc.}$$
We often omit tensor signs in writing elements of $T$
as is customary in noncommutive algebra, e.g.,
writing $x^3$ for $x\otimes x\otimes x$
and $xy$ for $x\otimes y$.

Suppose $P$ is a set of
filtered (nonhomogeneous) relations in degree 2,
$$P\subseteq T^0 \oplus T^1 \oplus T^2, $$
and let $I=(P)$ be the $2$-sided ideal in $T$ generated by $P$.
The quotient $A=T/I$ is a {\em nonhomogeneous quadratic
algebra}.
If $P$ consists of elements 
of homogeneous degree 2, i.e., $P\subseteq T^2$,
then $A$ is a {\em homogeneous quadratic algebra}.
Thus a quadratic algebra is just an algebra whose relations are generated
by (homogeneous or nonhomogenous) quadratic expressions.

We usually write each element of a finitely presented algebra
$A=T/I$ as a coset representative in $T$, suppressing
mention of the ideal $I$.  
Then a {\em $k$-basis} for  $A$ is a 
subset of $T$ representing cosets modulo $I$
which form a basis for $A$ as a $k$-vector space.
Some authors say a quadratic algebra has a {\em PBW basis}
if it has the same $k$-basis as a universal enveloping algebra,
i.e., if
$\{v_1^{a_1}\cdots v_n^{a_n}: a_i\in\NN\}$
is a basis for $A$ as a $k$-vector space.  Such algebras include
Weyl algebras, quantum/skew polynomial rings, some iterated Ore extensions,
some quantum groups, etc.

\subsection*{Priddy's PBW algebras}
Priddy~\cite{Priddy} gave a broader definition of PBW basis
for homogeneous quadratic algebras
in terms of any ordered basis of $V$
(say, $v_1< v_2 < \cdots < v_n$) 
in establishing the notion of Koszul algebras. 
(A quadratic algebra is {\em Koszul} if 
the boundary maps in its minimal free resolution have matrix entries
that  are linear forms;
see Section~\ref{Koszul}.)
Priddy first extended the ordering degree-lexicographically
to a monomial ordering on the tensor algebra  $T$,
where we regard pure tensors in $v_1, \ldots, v_n$ as monomials.
He then called a $k$-vector space basis for $A=T/I$
a {\em PBW basis} (and the algebra $A$ a {\em PBW algebra})
if the product of any two basis elements either lay
again in the basis or could be expressed modulo $I$
as a sum of larger elements in the basis.
In doing so, Priddy~\cite[Theorem 5.3]{Priddy} gave a class of Koszul algebras
which is easy to study:
%%%%%%%%%%%%%%%%%%%%%%%%
\begin{thm}
If a homogeneous quadratic algebra has a PBW basis, then it is Koszul.
\end{thm}
%%%%%%%%%%%%%%%%%%%%%%%%%%5

Polishchuk and Positselski
reframed Priddy's idea; we summarize their
approach (see~\cite[Chapter~4, Section~1]{PP})
using the  notion of leading monomial
$\LM$ of any element of $T$ written in terms of the basis
$v_1, \ldots, v_n$ of $V$.
Suppose the set of generating relations $P$ is a subspace of $T^2$.
Consider 
those monomials that are not divisible
by the leading monomial of any generating quadratic relation:
$$
{\mathcal{B}}_{P}=\{\text{monomials } m\in T:
\LM(a) \nmid m,\ \forall a\in P
\}\, .
$$
Polishchuk and Positselski
call $\mathcal{B}_P$ 
a {\em PBW basis} of the quadratic algebra $A$ (and $A$ a {\em PBW
algebra}) whenever $\mathcal{B}_P$ is a $k$-basis of $A$.

\subsection*{Gr\"obner bases}
Priddy's definition and the reformulation 
of Polishchuk and Positselski 
immediately call to mind the theory of
Gr\"obner bases.
Recall that a set $\GB$ of nonzero elements generating 
an ideal $I$ is called a (noncommutative) {\em Gr\"obner basis} if
the leading monomial of each nonzero element of $I$ is divisible
by the leading monomial of some element of $\GB$
with respect to a fixed monomial (i.e., term) ordering 
(see~\cite{Mora} or~\cite{Li2012}).
(Gr\"obner bases and Gr\"obner-Shirshov bases were developed
independently in various contexts by Shirshov~\cite{ShirshovOn62} in 1962, 
Hironaka~\cite{Hironaka} in 1964, Buchberger~\cite{BuchbergerThesis} in 1965,
Bokut'~\cite{Bokut} in 1976, 
and Bergman~\cite{Bergman} in 1978.)
A Gr\"obner basis $\GB$ is {\em quadratic}
if it consists of homogeneous elements of degree 2
(i.e., lies in $T^2$) and it is {\em minimal} if no proper subset is
also a Gr\"obner basis.  A version of the Composition-Diamond 
Lemma for associative algebras (see Section~\ref{diamond})
implies that if $\GB$ is a Gr\"obner basis for
$I$, then
$$\mathcal{B}_{\GB}=\{\text{monomials } m \in T: \,
\LM(a) \nmid m,\, \forall
a \in \GB\}
$$
is a $k$-basis for $A=T(V)/I$.

\vspace{2ex}
%%%%%%%%%%%%%%%%%%%
\begin{ex}
Let $A$ be the $\CC$-algebra generated by symbols $x,y$ with
a single generating relation $xy=y^2$.  
Set $V=\CC\text{-span}\{x,y\}$ and $P=\{xy-y^2\}$
so that $A=T(V)/(P)$.  
A Gr\"obner basis $\GB$ for the ideal
$I=(P)$ with respect to the degree-lexicographical
monomial ordering with $x<y$ is 
infinite:
$$\begin{aligned}
\GB&=\{yx^ny-x^{n+1}y: n\in \NN\},\\
\mathcal{B}_{P}&=\{\text{monomials } m\in T
\text{ that are not divisible by } y^2\},\\
\mathcal{B}_{\GB}&=\{\text{monomials } m\in T
\text{ that are not divisible by } yx^ny \text{ for any } 
n\in \NN\}.
\end{aligned}
$$
Hence, $A$ is not a PBW algebra using the ordering $x<y$
since $\mathcal{B}_{\GB}$ is a $\CC$-basis for
$A$ but $\mathcal{B}_{P}$ is not.

If we instead take some monomial ordering with $x>y$,
then $\GB=P$ is a Gr\"obner basis for the ideal $I=(P)$ 
and $\mathcal{B}_{\GB}=\mathcal{B}_{P}$ is a $\CC$-basis of $A$:
$$
\begin{aligned}
\mathcal{B}_{P}=\mathcal{B}_{\GB}&=
\{\text{monomials } m \in T
\text{ that are not divisible by } xy\}\\
&=\{y^ax^b:a,b\in \NN\}.
\end{aligned}
$$
Hence $A$ is a PBW algebra using the ordering $y<x$.
\end{ex}
%%%%%%%%%%%%%%%%%%%
\vspace{2ex}

\subsection*{Quadratic Gr\"obner bases}
How do the sets of monomials $\mathcal{B}_P$ and  $\mathcal{B}_{\GB}$
compare after fixing an appropriate monomial ordering?
Suppose $\GB$ is a minimal Gr\"obner basis for $I=(P)$ (which implies that
no element of $\GB$ has leading monomial dividing that of another).
Then $\mathcal{B}_{\GB}\subset \mathcal{B}_P$,
and the reverse inclusion holds whenever $\GB$ is quadratic
(since then $\GB$ must be a subset of the subspace $P$).
Since each graded piece of $A$ is finite dimensional over $k$, 
a PBW basis thus corresponds to a quadratic Gr\"obner basis:
$$\mathcal{B}_P \text{ is a PBW basis of }A
\iff
\mathcal{B}_{\GB}= \mathcal{B}_P
\iff
\GB \text{ is quadratic}.
$$
Thus authors sometimes call any algebra
defined by an ideal of relations with a
quadratic Gr\"obner basis a PBW algebra.
In any case (see~\cite{Anick},~\cite{BHV},~\cite{Froberg}):
\begin{thm}
Any quadratic algebra 
whose ideal of relations has a (noncommutative) quadratic
Gr\"obner basis is Koszul.
\end{thm}

Backelin (see~\cite[Chapter 4, Section 3]{PP}) gave an example
of a Koszul algebra defined by an ideal of relations with no
quadratic Gr\"obner basis.
Eisenbud, Reeves, and Totaro~\cite[p.~187]{ERT} 
gave an example
of a commutative Koszul algebra whose ideal of relations
does not have a quadratic Gr\"obner basis with respect to {\em any} ordering,
even after a change of basis (see also~\cite{Froberg}).

We relate Gr\"obner bases and PBW theorems for {\em nonhomogeneous}
algebras in Section~\ref{diamond}.

%%%%%%%%%%%%%%%%%%%%%%%%%%%%%%%%%%%%%%%%%%%%%%%%%%%%%%%%%%
%%%%%%%%%%%%%%%%%%%%%%%%%%%%%%%%%%%%%%%%%%%%%%%%%%%%%%%%%5
\section{Nonhomogeneous algebras: PBW deformations}\label{sec:nonhomdef}

Algebras defined by generators
and relations are not naturally graded, but merely filtered,
and often one wants to pass
to some graded or homogeneous version of the algebra for
quick information.  
There is more than one way to do this in general. 
The original PBW Theorem shows
that the universal enveloping algebra of a Lie algebra
has one natural homogeneous version.  Authors
apply this idea  to other algebras, saying that 
an algebra satisfies a {\em PBW property} 
when 
graded versions are isomorphic  
and call the original algebra a {\em PBW deformation}
of this graded version. We make these notions
precise in this section and relate them to the work of Braverman and
Gaitsgory and of Polishchuk and Positselski on Koszul algebras in the next section.

\subsection*{Filtered algebras}
Again, consider an algebra $A$ generated by a finite dimensional 
vector space $V$ over a field $k$ with some defining set of relations $P$. 
(More generally,
one might consider a module over a group algebra or some other $k$-algebra.)
Let $T=\bigoplus_{i\geq 0} T^i$ be the tensor algebra over $V$
and let $I=( P)$ be the two-sided ideal of relations
so that
$$A={T}/ {I} \, .$$
If $I$ is  homogeneous, then the quotient algebra
$A$ is graded.  In general, $I$ is
nonhomogeneous and the quotient algebra is only filtered, 
with $i$-th filtered component
$F^i(A)=F^i(T/I) = (F^i(T)+I)/I$ induced from the filtration
on $T$ obtained by assigning degree one to each vector in $V$
(i.e., $F^i(T)= T^0\oplus T^1 \oplus \ldots \oplus T^i$).

\subsection*{Homogeneous versions}
One associates to the filtered 
algebra $A$ two possibly different graded versions.
On one hand, we cross out lower order terms in the {\em generating}
set $P$ of relations to  obtain a 
homogeneous version of the original algebra.
On the other hand, we  
cross out lower order terms 
in each element of the {\em entire} ideal of relations.
Then {\em PBW conditions} 
are precisely those 
under which these two 
graded versions of the original algebra 
coincide, as we recall next.

The {\em associated graded
algebra} of $A$,
$$\gr (A) 
= \ \bigoplus_{i\geq 0} \quotient{{F}^i(A)}{{F}^{i-1}(A)} \, ,
$$
is a graded version of $A$ which
does not depend on the choice of generators $P$ of the ideal
of relations $I$.  (We set $F^{-1} = \{0\}$.)
The associated graded algebra may be realized concretely by
projecting each element in the ideal $I$ onto its
leading homogeneous part
(see Li~\cite[Theorem 3.2]{Li2012}):
$$
\gr \big(\quotient{T}{I}\big)
\ \cong\
\quotient{T}{ ( \LH(I))}\, ,
$$
where $\LH(S)=\{\LH(f): f \in S\}$ for any $S\subseteq T$
and $\LH(f)$ picks off the leading
(or highest) homogeneous part of $f$ in the graded algebra $T$. 
(Formally, $\LH(f)=f_d$ for $f=\sum_{i=1}^d f_i$ with 
each $f_i$ in $T^i$ and $f_d$ nonzero.) 
Those looking for a shortcut may be tempted
instead simply to project elements of the generating set
$P$ onto their leading homogeneous parts.
A natural surjection (of graded algebras) always arises from this
homogeneous version of $A$ determined by $P$
to the associated graded algebra of $A$:
$$\quotient{T}{( \LH(P)) }
\twoheadrightarrow 
\gr \big(\quotient{T}{I}\big)\ .$$ 

\subsection*{PBW deformations}
We say the algebra $T/I$
is a {\em PBW deformation} of its homogeneous version
$T/ ( \LH(P))$
(or satisfies the {\em PBW property}
with respect to $P$) 
when the above surjection is also injective,
i.e., when the associated graded algebra
and the homogeneous algebra determined by $P$
coincide (see~\cite{BG}):
$$
\quotient{T}{( \LH(I) ) }
\cong
\gr \Big(\quotient{T}{I}\Big)\ 
\cong
\quotient{T}{( \LH(P))} \,
.$$
In the next section, we explain the connections among 
PBW deformations, graded (and formal) deformations, and Hochschild cohomology.

In this language, 
the original PBW Theorem for universal enveloping algebras
asserts that the set
$$
P=\{v\ot w-w\ot v-[v,w]:v,w\in V\}
$$
gives rise to a quotient algebra $T/(P)$ that is a PBW deformation
of the commutative polynomial ring $S(V)$,
for $V$ the underlying vector space of a Lie algebra.
Here, each element of $V$ has degree $1$ so that the relations
are nonhomogeneous of degree $2$ and $T/(P)$ is a nonhomogenous quadratic
algebra.

We include an example next to show
how the PBW property depends on choice of generating relations
$P$ defining the algebra $T/I$.  
(But note that if $A$ satisfies the PBW property with respect
to some generating set $P$ of relations, then the subspace $P$
generates is unique; 
see~\cite[Proposition 2.1]{PBWQuadratic}.)

\vspace{2ex}
%%%%%%%%%%%%%%%%%%%%%%%%%%%%%%%%%%%%%%%%%%%%%%%55
\begin{ex}\label{cuteexample}{
We mention a filtered algebra that exhibits
the PBW property
with respect to one generating set of relations but not another.
Consider the (noncommutative) algebra $A$ generated by symbols $x$ and $y$
with defining relations $xy=x$ and $yx=y$:
$$ A=\quotient{k\langle x,y \rangle }{(xy-x, yx-y)}\, ,$$
where $k\langle x,y\rangle$ is the free $k$-algebra generated
by $x$ and $y$.
The algebra $A$ {\em does  not} satisfy the PBW property with
respect to the generating relations $xy-x$ and $yx-y$.
Indeed, the relations imply that $x^2=x$ and $y^2=y$
in $A$ and thus the associated graded algebra $\gr(A)$ is trivial
in degree two while 
the homogeneous version of $A$ is not (as $x^2$ and $y^2$ represent
nonzero classes).
The algebra $A$
{\em does} exhibit the PBW property with respect to the larger generating set
$\{xy-x, yx-y, x^2-x, y^2-y\}$ since 
$$\gr A \cong 
\quotient{k\langle x,y\rangle}{(xy, yx, x^2, y^2)}\, . $$
Examples~\ref{cuteexampleGroebner} and~\ref{cuteexampleDiamond}
explain this recovery of the PBW property in terms of Gr\"obner bases
and the Composition-Diamond Lemma.}\end{ex}
%%%%%%%%%%%%%%%%%%%%%%%%%%%%%%%%%%%%%%%%%%%%%%%%%%%%%%%%
%%%%%%%%%%%%%%%%%%%%%%%%%%%%%%%%%%%%%%%%%%%%%%%%%%%%%%%%%%
%%%%%%%%%%%%%%%%%%%%%%%%%%%%%%%%%%%%%%%%%%%%%%%%%%%%%%%%%%
%%%%%%%%%%%%%%%%%%%%%%%%%%%%%%%%%%%%%%%%%%%%%%%%%%%%%%%%%%
\section{Deformation Theory and Hochschild cohomology}\label{defHH}

In the last section, we saw that an algebra defined by
nonhomogeneous relations is called a {\em PBW deformation}
when the homogeneous version determined by generating relations
coincides with its associated graded algebra.
How may one view formally the original nonhomogeneous algebra
as a {\em deformation} of its homogeneous version?
In this section, we begin to 
fit PBW deformations into the theory of algebraic deformations.
We recall the theory of deformations
of algebras and Hochschild cohomology, a homological tool used
to predict deformations and prove PBW properties.

%%%%%%%%%%%%%%%%%%%%%%%%%%%%%%%%%%%%%%%%
\subsection*{Graded deformations}

Let $t$ be a formal parameter.  
A {\em graded deformation} of a graded $k$-algebra $A$ 
is a graded associative $k[t]$-algebra $A_t$ (for $t$ in degree 1)
which is isomorphic to $A[t]=A\ot_k k[t]$ as a $k[t]$-module
with $$A_t|_{t=0} \cong A .$$
If we specialize $t$ to an element of $k$ in the algebra $A_t$, 
then we may no longer
have a graded algebra, but a filtered algebra instead.

PBW deformations may be viewed as graded deformations:
Each PBW deformation is a graded deformation of its homogeneous version
with parameter $t$ specialized to some element of $k$.
Indeed, given a finitely generated algebra $A=T/(P)$,
we may insert a formal parameter $t$ of degree 1
throughout the defining relations $P$ to make each relation homogeneous
and extend scalars to $k[t]$;
the result yields a graded algebra $B_t$ over $k[t]$
with 
$A=B_t|_{t=1}$ and $B= B_t|_{t=0}$, the homogeneous version of $A$.  
One may verify that if $A$ satisfies the PBW property, then
this interpolating algebra $B_t$ also satisfies a PBW condition over $k[t]$
and that $B_t$ and $B[t]$ are isomorphic as $k[t]$-modules.
Thus as $B_t$ is an associative graded algebra,
it defines a graded deformation of $B$.

Suppose $A_t$ is a graded deformation
of a graded $k$-algebra $A$.
Then up to isomorphism, $A_t$ is just the vector space $A[t]$ together
with some associative multiplication given by 
\begin{equation}\label{associativemultiplication}
  a * b = ab + \mu_1(a\ot b)t + \mu_2(a\ot b) t^2+\cdots,
\end{equation}
where $ab$ is the product of $a$ and $b$ in $A$ and for each $i$,
and each $\mu_i$ is a linear map from $A\otimes A$ to $A$ of degree $-i$,
extended to be $k[t]$-linear.
The degree condition on the maps $\mu_i$ are forced by the fact
that $A_t$ is graded for $t$ in degree 1.
(One sometimes considers a {\em formal} deformation, defined
over formal power series $k[[t]]$ instead of polynomials $k[t]$.) 

The condition that the  multiplication $*$ in $A[t]$ be associative imposes
conditions on the functions $\mu_i$ which are often
expressed using Hochschild cohomology. For example, comparing
coefficients of $t$ in the equation $(a*b)*c = a*(b*c)$, we see that 
$\mu_1$ must satisfy
\begin{equation}\label{cocyclecondition}
   a\mu_1(b\ot c) + \mu_1(a\ot bc) = \mu_1(ab\ot c) + \mu_1(a\ot b)c
\end{equation}
for all $a,b,c\in A$. We see below that this condition
implies that $\mu_1$ is a Hochschild 2-cocycle.
Comparing coefficients of $t^2$ yields
a condition on $\mu_1,\mu_2$ called the 
{\em first obstruction}, comparing coefficients 
of $t^3$ yields a condition on $\mu_1, \mu_2,\mu_3$
called the {\em second obstruction}, and so on.
(See~\cite{Gerstenhaber}.) 

\subsection*{Hochschild cohomology}
Hochschild cohomology is a generalization of group cohomology
well suited to noncommutative algebras.  It 
gives information about an algebra $A$ viewed as a bimodule over itself,
thus capturing right and left multiplication,
and predicts possible multiplication maps $\mu_i$ that could
be used to define a deformation of $A$.
One may define the Hochschild cohomology of a $k$-algebra
concretely as Hochschild cocycles modulo Hochschild coboundaries
by setting 
$$
\text{Hochschild $i$-cochains} =
\{ \text{linear functions } \phi: 
\underbrace{A\otimes \cdots \otimes A}_{i-\text{times}}
\rightarrow A \}
$$
(i.e., multilinear maps $A\times \cdots\times A \rightarrow A$)
with linear boundary operator 
$$\delta_{i+1}^*: i\text{-cochains}
\rightarrow (i+1)\text{-cochains}$$
given by
$$
\begin{aligned}
(\delta_{i+1}^*\phi)&(a_0\otimes \cdots\otimes a_{i})
=\\
& a_0\phi(a_1\otimes \cdots\otimes a_{i}) 
+\sum_{0\leq j\leq i-1} (-1)^{j+1} 
\phi(a_0\otimes\cdots\otimes a_{j-1}\otimes a_j a_{j+1}\otimes a_{j+2}
\otimes\cdots\otimes a_{i})\\
&\hphantom{xxxxxxxxxxxxxxxxxxxxxxxxxxxxxx}
 + (-1)^{i+1}\phi(a_0\otimes\cdots\otimes a_{i-1}) a_{i}\, .
\end{aligned}
$$
We identify $A$ with $\{0\text{-cochains}\}$.
Then
$$
\HH^i(A):=\ker \delta_{i+1}^* / \text{Im } \delta_i^*\ .
$$

We are interested in other concrete realizations of Hochschild cohomology
giving isomorphic cohomology groups.
Formally, we view any $k$-algebra $A$ as a bimodule over itself,
i.e., a right $A^e$-module where $A^e$ is its 
enveloping algebra, $A\ot A^{op}$, for 
$A^{\op}$ the 
opposite algebra of $A$. The Hochschild cohomology of $A$ is
then just
$$
\HH^{\DOT}(A)=\text{Ext}^{\DOT}_{A^e}(A,A).
$$
This cohomology is often computed
using the $A$-bimodule bar resolution of $A$:
\begin{equation}\label{res-bar}
  \cdots \stackrel{}{\longrightarrow} 
   A^{\ot 4}\stackrel{\delta_2}{\longrightarrow} A^{\ot 3}
  \stackrel{\delta_1}{\longrightarrow} A^{\ot 2}
  \stackrel{\delta_0}{\longrightarrow} A
   \longrightarrow 0  ,
\end{equation}
where $\delta_0$ is the multiplication in $A$, and, for each $i\geq 1$, 
$$
   \delta_i(a_0\ot\cdots\ot a_{i+1}) =
     \sum_{j=0}^{i} (-1)^j a_0\ot \cdots\ot a_j a_{j+1}\ot \cdots\ot a_{i+1} 
$$
for $a_0,\ldots, a_{i+1}$ in $A$.  We take 
the homology of this complex after dropping the initial term $A$ 
and applying $\Hom_{A\ot A^{\text{op}}}(-,A)$ to obtain the above
description of Hochschild cohomology
in terms of Hochschild cocycles and coboundaries,
using the identification
$$
\Hom_{A\ot A^{\text{op}}}(A\ot A^{\otimes i}\ot A, A)
\cong
\Hom_{k}(A^{\otimes i}, A).
$$ 

%%%%%%%%%%%%%%%%%%%%%%%%%%%%%%%%%%%%%%%%%%%%%%%%%%%%%%%%%%%%
%%%%%%%%%%%%%%%%%%%%%%%%%%%%%%%%%%%%%%%%%%%%%%%%%%%%%%%%%%%%%%%%5
\section{Koszul algebras}\label{Koszul}

We wish to extend the original PBW Theorem for
universal enveloping algebras to other
nonhomogeneous quadratic algebras.
When is a given algebra
a PBW deformation of another well-understood and
well-behaved algebra?  
Can we replace the polynomial algebra in 
the original PBW theorem by any homogeneous quadratic algebra, provided it is 
well-behaved in some way?
We turn to Koszul algebras as a wide class of quadratic algebras
generalizing the class of polynomial algebras.
In this section, we briefly recall the definition of a Koszul
algebra.  

\subsection{Koszul complex}
The algebra $S$ is a {\em Koszul algebra} if the underlying field $k$ admits
a linear $S$-free resolution, i.e., 
one with boundary maps given by matrices whose 
entries are linear forms.  Equivalently, $S$ is a Koszul algebra
if the following  complex of left
$S$-modules is acyclic:
\begin{equation}\label{res-koszul}
  \cdots\longrightarrow K_3(S)\longrightarrow K_2(S)\longrightarrow
   K_1(S)\longrightarrow K_0(S)\longrightarrow k\longrightarrow 0
\end{equation}
where $K_0(S) = S$, $K_1(S)=S\ot V$, $K_2(S)=S\ot R$, and 
for $i\geq 3$, 
$$
   K_i(S) = S\ot \left(\,\,\bigcap_{j=0}^{i-2} V^{\ot j}\ot R\ot
     V^{\ot (i-2-j)} \right) .
$$
The differential is that inherited from the bar resolution
of $k$ as an $S$-bimodule,
\begin{equation}\label{res-bar-k}
  \cdots \stackrel{\del_4}{\longrightarrow} 
   S^{\ot 4}\stackrel{\del_3}{\longrightarrow} S^{\ot 3}
  \stackrel{\del_2}{\longrightarrow} S^{\ot 2}
  \stackrel{\del_1}{\longrightarrow} S
   \stackrel{\varepsilon}{\longrightarrow} k \longrightarrow 0  ,
\end{equation}
where $\varepsilon$ is the augmentation ($\varepsilon(v)=0$ for all $v$ in $V$)
and for each $i\geq 1$, 
$$
   \del_i(s_0\ot\cdots\ot s_{i}) =
      (-1)^i \varepsilon(s_i) s_0\ot \cdots \ot s_{i-1} +
     \sum_{j=0}^{i-1} (-1)^j s_0\ot \cdots\ot s_j s_{j+1}\ot \cdots\ot s_{i} .
$$
(Note that for each $i$, $K_i(S)$ is an $S$-submodule of $S^{\ot (i+1)}$.)

\subsection*{Bimodule Koszul complex}
Braverman and Gaitsgory gave an equivalent definition of Koszul algebra
via  the bimodule 
Koszul complex: Let 
\begin{equation}\label{K-tilde}
\widetilde{K}_i(S) = K_i(S)\ot S ,
\end{equation}  
an $S^e$-module (equivalently $S$-bimodule) where $S^e=S\ot S^{op}$.
Then $\widetilde{K}_{\DOT}(S)$ embeds into the bimodule bar resolution
(\ref{res-bar}) whose $i$-th term is $S^{\ot (i+2)}$, and 
$S$ is Koszul if and only if $\widetilde{K}_{\DOT}(S)$ is a
bimodule resolution of $S$.
Thus we may obtain the 
Hochschild cohomology $\HH^{\DOT}(S)$ of $S$
(which contains information about its deformations)
by applying $\Hom_{S^e}( - , S)$ either to the Koszul resolution
$\widetilde{K}_{\DOT}(S)$ or
to the bar resolution (\ref{res-bar})
of $S$ as an $S^e$-module (after dropping the initial nonzero terms of each)
and taking homology.
We see in the next section how these 
resolutions and the resulting cohomology are used in homological
proofs of a  generalization of the PBW Theorem
from~\cite{BG,PP,Positselski}.

%%%%%%%%%%%%%%%%%%%%%%%%%%%%%%%%%%%%%%%%%%%%%%%%%%%%%%%%%%
%%%%%%%%%%%%%%%%%%%%%%%%%%%%%%%%%%%%%%%%%%%%%%%%%%%%%%%%%%
%%%%%%%%%%%%%%%%%%%%%%%%%%%%%%%%%%%%%%%%%%%%%%%%%%%%%%%%%%
\section{Homological methods and  deformations of Koszul algebras}
\label{BG}

Polishchuk and Positselski~\cite{PP,Positselski}
and Braverman and Gaitsgory~\cite{BG} 
extended the idea of the original PBW Theorem for
universal enveloping algebras to other
nonhomogeneous quadratic algebras by replacing the polynomial algebra in 
the theorem by an arbitrary Koszul algebra.
They stated conditions for
a version of the original PBW Theorem to hold
in this greater generality and gave
homological proofs. 
(Polishchuk and Positselski~\cite{PP} in fact gave two proofs, one homological
that goes back to Positselski~\cite{Positselski}
and another using distributive lattices.)
We briefly summarize these two homological approaches in this section
and discuss generalizations.

%%%%%%%%%%%%%%%%%%%%%%%%%%%%%%%%
\subsection*{Theorem of Polishchuk and Positselski,
Braverman and Gaitsgory}
As in the last sections, 
let $V$ be a finite dimensional vector space over a field $k$
and let $T$ be its tensor algebra over $k$
with $i$-th filtered component $F^i(T)$.
Consider a subspace $P$ of $F^2(T)$ defining a nonhomogeneous
quadratic algebra $$A=\quotient{T}{(P)}\ .$$
Let $R=\LH(P)\cap T^2$ be the projection of $P$ onto the
homogeneous component of degree 2, and 
set  $$S=\quotient{T}{(R)},$$
a homogeneous quadratic algebra (the homogeneous version
of $A$ as in Section~\ref{sec:nonhomdef}).
Then $A$ is a PBW deformation of $S$ when $\gr A$
and $S$ are isomorphic as graded algebras.

Braverman and Gaitsgory and also Polishchuk and Positselski 
gave a generalization of the PBW Theorem~\cite{BG,PP,Positselski} as follows:

%%%%%%%%%%%
\begin{thm}\label{IJ}
Let $A$ be a nonhomogeneous quadratic algebra, $A = T/(P)$, and $S=T/(R)$
its corresponding homogeneous quadratic algebra. 
Suppose $S$ is a Koszul algebra. 
Then $A$ is a PBW deformation of $S$ if, and only if, the following
two conditions hold:

(I) $ \ P\cap F^1(T) = \{0\}$, and

(J) $ \ (F^1(T)\cdot P \cdot F^1(T))\cap F^2(T) = P$.
\end{thm} 

We have chosen the notation of Braverman and Gaitsgory.
The necessity of  conditions (I) and (J) can be seen by direct algebraic manipulations.
Similarly,  direct computation 
shows  that if (I) holds, then (J) is equivalent 
to (i), (ii), and (iii) of Theorem~\ref{thm:BG}  below. 
Braverman and Gaitsgory used 
algebraic deformation theory to show that these conditions are also sufficient. 
Polishchuk and Positselski used properties of an explicit complex defined 
using the Koszul dual of $S$. 
The conditions (i), (ii), (iii) facilitate these connections to 
homological algebra, and 
they are easier in practice to verify than checking (J) directly.
But in order to state these conditions, we require a
canonical decomposition for  elements of $P$:
Condition~(I) of Theorem~\ref{IJ} implies that every element of $P$
can be written as the sum of a nonzero element of $R$ (of degree $2$), 
a linear term, and a constant term,
i.e., there exist linear functions
$\alpha: R\rightarrow V$, $\beta: R\rightarrow k$ for which 
$$P=\{r-\alpha(r)-\beta(r)\mid r\in R\}.$$
One may then rewrite 
Condition~(J) and
reformulate Theorem~\ref{IJ} as follows.

%%%%%%%%%%%%%%%%%%%%%%%%%%%%%%%%%%
\begin{thm}\label{thm:BG}
Let $A$ be a nonhomogeneous quadratic algebra, $A=T/(P)$, and $S=T/(R)$
its corresponding homogeneous quadratic algebra.
Suppose $S$ is a Koszul algebra.
Then $A$ is a PBW deformation of $S$ if, and only if, the following
conditions hold:

(I)  $ \ \ \ P\cap F^1(T) = \{0\}$,

(i) $\ \ \ \Ima (\alpha\ot \id - \id\ot \alpha) \subseteq R$,

(ii) $ \ \ \alpha\circ (\alpha\ot\id - \id\ot \alpha) =
      - (\beta\ot\id -\id\ot \beta)$,

(iii) $ \ \beta\circ (\alpha\ot\id - \id\ot \alpha) = 0$,

\noindent
where the maps $\alpha\ot\id -\id\ot\alpha$ and $\beta\ot\id - \id\ot\beta$ are
defined on the subspace $(R\ot V)\cap (V\ot R)$ of $T$. 
\end{thm}
%%%%%%%%%%%%%%%%%%%%%%%%%%%%%%%%%%%%%%%%%%%%%%

We explain next how the original PBW Theorem is a consequence
of Theorem~\ref{thm:BG}. 
Indeed, Polishchuk and 
Positselski~\cite[Chapter 5, Sections 1 and 2]{PP}
described the ``self-consistency conditions'' (i), (ii), and (iii)
of the theorem
as  generalizing the Jacobi identity
for Lie brackets.

\vspace{2ex}
%%%%%%%%%%%%%%%%%%%%%%%%%%%
\begin{ex}
Let $\mathfrak g$ be a finite dimensional complex Lie algebra, $A=U({\mathfrak{g}})$
its universal enveloping algebra, and $S=S({\mathfrak{g}})$.
Then $R$ has $\CC$-basis all $v\ot w - w\ot v$ for $v,w$ in $V$,
 and $\alpha(v\ot w-w\ot v)
= [v,w]$, $\ \beta\equiv 0$.
Condition (I) is equivalent to antisymmetry of the bracket.
Condition (J) is equivalent to the Jacobi identity, with (i), (ii) expressing
the condition separately in each degree in the tensor algebra
($\beta \equiv 0$ in this case). More generally, 
there are examples with $\beta\not\equiv 0$, 
for instance, the Sridharan enveloping algebras~\cite{Sridharan}.
\end{ex} 
%%%%%%%%%%%%%%%%%%%%%%%%%%
\vspace{2ex}

%%%%%%%%%%%%%%%%%%%%%%%%%%%%%%%%
\subsection*{Homological proofs}
We now explain how Braverman and Gaitsgory 
and Polishchuk and Positselski 
used  algebraic deformation
theory and Hochschild cohomology
to prove that the conditions of Theorem~\ref{thm:BG} are 
sufficient.
Braverman and Gaitsgory constructed
a graded deformation $S_t$ interpolating between $S$ and $A$
(i.e., with $S=S_t|_{t=0}$ and $A=S_t|_{t=1}$),
implying that $\gr (A)\cong S$ as graded algebras.
They constructed the deformation $S_t$ as follows.
\begin{itemize}
\item
They identified $\alpha$
with a choice of first multiplication map $\mu_1$ and $\beta$
with a choice of second multiplication map $\mu_2$, via the canonical
embedding of the bimodule Koszul resolution (\ref{K-tilde}) into the
bar resolution (\ref{res-bar}) of $S$. 
(In order to do this, one must extend $\alpha,\beta$ 
(respectively, $\mu_1,\mu_2$) to
be maps on a larger space via an isomorphism 
$\Hom_k(R,S)\cong \Hom_{S^e}(S\ot R\ot S,S)$
(respectively, $\Hom_k(S\ot S,S)\cong \Hom_{S^e}(S^{\ot 4},S)$.)
\item 
Condition (i) is then seen to be equivalent to $\mu_1$ being a Hochschild
2-cocycle (i.e., satisfies Equation~(\ref{cocyclecondition})). 
\item
Condition (ii) is equivalent to the vanishing of the first
obstruction. 
\item
Condition (iii) is equivalent to the vanishing of
the second obstruction. 
\item
All other obstructions vanish automatically for a 
Koszul algebra due to the structure of its Hochschild cohomology 
(see~\cite{BG}).  
\item
Thus there exist maps $\mu_i$ for $i>2$
defining an associative multiplication $*$ 
(as in Equation~(\ref{associativemultiplication}))
on $S[t]$.
\end{itemize}

Positselski~\cite[Theorem~3.3]{Positselski} (see also \cite[Proposition~5.7.2]{PP})
 gave a different homological proof of Theorem~\ref{thm:BG}.
Let $B$ be the Koszul dual $S^{!}:= \Ext^*_S(k,k)$ of $S$.
Then $S\cong B^{!}:= \Ext^*_B(k,k)$. Polishchuk defined a complex
whose terms are the same as those in the bar resolution of $B$ but with
boundary maps
modified using the functions $\alpha: R\rightarrow V$,
$\beta: R\rightarrow k$ by first identifying $\beta$ with an
element $h$ of $B^2$ and $\alpha$ with a 
dual to a derivation $d$ on $B$.
The conditions (i), (ii), and (iii) on $\alpha,\beta$ correspond to
conditions on $d,h$, under which Positselski called $B$ a CDG-algebra.
The idea is that CDG-algebra structures on $B$ are dual to PBW deformations
of $S$. Positselski's proof relies on the Koszul property of $S$
(equivalently of $B$) to imply collapsing of a spectral sequence
with $E^2_{p,q} = \Ext^{-q,p}_B(k,k)$. The sequence converges
to the homology of the original complex for $B$. Koszulness
implies the only nonzero terms occur when $p+q=0$, and we are
left with the homology of the total complex in degree 0. By its definition
this is simply the algebra $A$, and it follows that $\gr A\cong B^{!}\cong S$. 

%%%%%%%%%%%%%%%%%%%%%%%%%%%%%%%%
\subsection*{Generalizations and extensions}
Theorem~\ref{thm:BG} describes nonhomogeneous quadratic algebras
whose quadratic versions are Koszul.  What if one replaces the
underlying field by an arbitrary ring? 
Etingof and Ginzburg~\cite{EtingofGinzburg} noted that Braverman and Gaitsgory's
proof of Theorem~\ref{thm:BG} is in fact  valid more generally for 
Koszul rings over semisimple subrings 
as defined by Beilinson, Ginzburg, and Soergel~\cite{BGS}.
They chose their semisimple subring to be the complex group algebra $\CC G$ of a finite group $G$ acting symplectically and their Koszul ring to 
be a polynomial algebra $S(V)$. They 
were interested  in the case $\alpha\equiv 0$ for their applications
to symplectic reflection algebras (outlined in Section~\ref{SRA} below).
Halbout, Oudom, and Tang~\cite{HOT} state a 
generalization of Theorem~\ref{thm:BG} in this setting that allows 
nonzero $\alpha$ (i.e., allows relations defining the algebra $A$ 
to set commutators of vectors in $V$ to a combination of group algebra
elements and vectors).
A proof using the Koszul ring theory of Beilinson, Ginzburg, and 
Soergel and the results of Braverman and Gaitsgory
is outlined in our paper~\cite{doa} 
for arbitrary group algebras over
the complex numbers.
We also included a second proof there for group algebras
over arbitrary fields (of characteristic not 2) using
the Composition-Diamond Lemma (described in the next section),
which has the advantage that it is characteristic free. 
We adapted the program of Braverman and Gaitsgory 
to arbitrary
nonhomogeneous quadratic algebras and Koszul rings
defined over non-semisimple rings in~\cite{PBWQuadratic}, 
including group rings $kG$ where the characteristic of
$k$ divides the order of the group $G$.

The theory of Braverman and Gaitsgory was further generalized to algebras
that are $N$-Koszul (all relations homogeneous of degree $N$ plus a 
homological condition) over semisimple or von Neumann regular rings
by a number of authors
(see~\cite{BergerGinzburg,FloystadVatne,HSS}).
Cassidy and Shelton~\cite{CassidyShelton} generalized the theory of Braverman and
Gaitsgory  in a different direction, to graded algebras over a
field satisfying a particular homological finiteness condition
(not necessarily having all relations in a single fixed degree). 

%%%%%%%%%%%%%%%%%%%%%%%%%%%%%%%%%%%%%%%%%%%%%%%%%%%%%%
\section{The Composition-Diamond 
Lemma and Gr\"obner basis theory}\label{diamond}

PBW theorems are often proven using diamond or composition lemmas
and the theory of (noncommutative) Gr\"obner bases.
Diamond lemmas predict existence of
a canonical normal form in a mathematical system.
Often one is presented with various
ways of simplifying an element to obtain a normal form.
If two different ways of rewriting the original element
result in the same desired reduced expression, one is reminded
of diverging paths meeting 
like the sides of the shape of a diamond.
Diamond lemmas often originate from Newman's Lemma~\cite{Newman} 
for graph theory. 
Shirshov (see~\cite{ShirshovOn62}
and~\cite{ShirshovSome62})
gave a general version
for Lie polynomials in 1962 which 
Bokut' (see~\cite{Bokut} and~\cite{BokutChen}) extended to 
associative algebras in 1976,
using the term ``Composition Lemma.''
Around the same time (Bokut' cites a preprint by Bergman),
Bergman~\cite{Bergman} developed a similar result which he 
instead called the Diamond Lemma.

Both the Diamond Lemma
and Composition Lemma are easy to explain 
but difficult to state precisely without the formalism
absorbed by Gr\"obner basis theory.
In fact,
the level of rigor
necessary to carefully state and prove these results
can be the subject of debate.
Bergman himself writes that the lemma
``has been considered obvious and used freely by some
ring-theorists... but others seem unaware
of it and write out
tortuous verifications.''
(Some authors are reminded
of life in a lunatic asylum 
(see~\cite{HellstromSilvestrov}) when making the basic
idea rigorous.)  We leave careful 
definitions to any one of numerous texts 
(for example, see~\cite{BMM} or~\cite{AlgorithmicMethods}) 
and instead present the intuitive idea behind
the result developed by Shirshov, Bokut', and Bergman.
%%%%%%%%%%%%%%%%%%%%%%%%%%%%%%%%%%%%%%%%%%%%%%%
\subsection*{The Result of Bokut' (and Shirshov)}
We first give the original result of Bokut' (see~\cite[Proposition 1 and Corollary 1]{Bokut}),
who used a degree-lexicographical monomial ordering
(also see~\cite{BokutKukin}).

%%%%%%%%%%%%%%%%%%%%%%%%%%%%%%%%%%%%%%
\begin{namedthm}[Original Composition Lemma]
Suppose a set of relations $P$ defining a $k$-algebra
$A$ is ``closed under composition.''
Then the set of monomials that do not contain
the leading monomial of any element of $P$ as a subword
is a $k$-basis of $A$. 
\end{namedthm}
%%%%%%%%%%%%%%%%%%%%%%%%%%%%%%%%%%%%%%%%%%%%%%%%%

Before explaining the notion of ``closed under composition,''
we rephrase the results of
Bokut' in modern language using Gr\"obner bases 
to give a PBW-like basis as in 
Section~\ref{homogeneous} 
(see~\cite{Green94}, or~\cite{Mora}, for example).
Fix a monomial ordering
on a free $k$-algebra $T$ and again write $\LM(p)$ for the leading
monomial of any $p$ in $T$.
We include the converse of the lemma
which can be deduced from the work 
of Shirshov and Bokut' and was given explicitly by Bergman,
who used
arbitrary monomial orderings.

%%%%%%%%%%%%%%%%%%%%%%%%%%%%%%
\begin{namedthm}[Gr\"obner basis version of Composition Lemma]
The set $P$ is a (noncommutative) Gr\"obner basis
of the ideal $I$ it generates if and only if
$$\mathcal{B}_{P}=\{\text{monomials } m \text{ in } T: \
m \text{ not divisible by any }
\LM(p),\ p \in P\}
$$
is a $k$-basis for the algebra $A=T/I$.
\end{namedthm}
%%%%%%%%%%%%%%%%%%%%%%%%%%%

\vspace{2ex}
%%%%%%%%%%%%%%%%%%%%%%%%%%%%%%%%%%%%%%%%%%55
\begin{ex}\label{cuteexampleGroebner}
Let $A$ be the $k$-algebra generated by symbols
$x$ and $y$ and relations $xy=x$ and $yx=y$ (Example~\ref{cuteexample}):
$$A=\quotient{k\langle x, y\rangle}{(xy-x, yx-y)}\, .$$
Let $P$ be the set of defining relations, $P=\{xy-x, yx-y\}$,
and consider the degree-lexicographical monomial ordering with $x>y$.
Then $P$ is {\em not} a Gr\"obner basis of the ideal it generates
since
$x^2-x=x(yx-y)-(xy-x)(x-1)$ lies in the ideal $(P)$ and has
leading monomial $x^2$, which does not lie in the ideal
generated by the leading monomials of the elements of $P$.
Indeed, $\mathcal{B}_{P}$ contains both $x^2$ and $x$
and hence can not be a basis for $A$.
We set $P'=\{xy-x, yx-y, x^2-x, y^2-y\}$ to obtain a Gr\"obner basis of $(P)$.
Then 
$\mathcal{B}_{P'}=\{\text{monomials } m : 
m \text{ not divisible by } xy, yx, x^2, y^2 \}$
is a $k$-basis for the algebra $A$.
\end{ex}
%%%%%%%%%%%%%%%%%%%%%%%%%%%%%%%%%%%%%
\vspace{2ex}

%%%%%%%%%%%%%%%%%%%%%%%%%%%%%%%%%%%%%%%%%%
\subsection*{Resolving ambiguities}
Bergman focused on the problem of
resolving ambiguities that arise when trying to rewrite
elements of an algebra using different defining relations.
Consider a $k$-algebra $A$ defined by a finite set of generators 
and a finite set of relations 
$$
m_1=f_1,\ m_2=f_2,\ \ldots,\ m_k=f_k\, ,
$$
where the $m_i$ are monomials (in the set of generators of $A$)
and the $f_i$ are linear combinations of monomials.
Suppose we prefer the right side of our relations and
try to eradicate the $m_i$ whenever possible in writing
the elements of $A$ in terms of its generators.
Can we define the notion of a canonical form for every element 
of $A$ by deciding to replace each $m_i$ by $f_i$ whenever possible?
We say an expression for an element of $A$ is {\em reduced}
if no $m_i$ appears (as a subword anywhere), i.e.,
when no further replacements using the defining relations of $A$
are possible.
The idea of a {\em canonical form} for $A$ then makes sense
if the set of reduced expressions is a $k$-basis for $A$,
i.e.,
if every element can be written uniquely in reduced form.

A natural ambiguity arises: If a monomial $m$ contains both $m_1$ and $m_2$
as (overlapping) subwords, do we ``reduce'' first $m_1$ to $f_1$ or
rather first $m_2$ to $f_2$ by replacing?  (In the last example,
the word $xyx$ contains overlapping subwords $xy$ and $yx$.)
If the order of application
of the two relations does not matter and we end up with the same reduced
expression, then we say the
{\em (overlap) ambiguity was resolvable}.
The Composition-Diamond Lemma states that 
{\em knowing certain ambiguities resolve is enough
to conclude that a canonical normal form
exists for all elements in the algebra}.

\vspace{2ex}
%%%%%%%%%%%%%%%%%%%%%%%%%%%%%%%%%%%%%%%%%%55
\begin{ex}\label{cuteexampleDiamond}
Again, let $A$ be the $k$-algebra generated by symbols
$x$ and $y$ and relations $xy=x$ and $yx=y$ (Example~\ref{cuteexample}). 
We decide to 
eradicate $xy$ and $yx$ whenever possible in expressions
for the elements of $A$ using just the defining relations.
On one hand, we may reduce $xyx$ to $x^2$ (using the first relation);
on the other hand, we may reduce $xyx$ to $xy$ (using
the second relation) then to $x$ (using the first relation).
The words $x$ and $x^2$ can not be reduced further
using just the defining relations, so we consider them both ``reduced''.
Yet they represent
the same element $xyx$ of $A$.
Thus, a canonical ``reduced'' form does not make sense given this choice
of defining relations for the algebra $A$.
\end{ex}
%%%%%%%%%%%%%%%%%%%%%%%%%%%%%%%%%%%%%
\vspace{2ex}

%%%%%%%%%%%%%%%%%%%%%%%%%%%%%%%%%%%%%%%%%%
\subsection*{The result of Bergman}
One makes the above notions precise by introducing a monomial ordering
and giving formal definitions for ambiguities,
reduction, rewriting procedures, resolving, etc.
We consider the quotient algebra $A=T/(P)$ where
$T$ (a tensor algebra) is the free $k$-algebra on the generators of $A$
and $P$ is a (say) finite set of generating relations.
We single out a monomial $m_i$ in each generating relation, writing
$$P=\{ m_i-f_i: 1\leq i\leq k\}\, ,$$
and choose a monomial ordering so that $m_i$ is the leading monomial
of each $m_i-f_i$ (assuming such an ordering exists).
Then the reduced elements are exactly those spanned by 
$\mathcal{B}_{P}$.  If all the ambiguities among elements of
$P$ are resolvable, we obtain a PBW-like basis, but Bokut' and Bergman
give a condition that is easier to check.
Instead of choosing to replace monomial $m_1$ by $f_1$ or monomial $m_2$ by $f_2$
when they both appear as subwords of a monomial $m$,
we make {\em both} replacements 
separately and take the difference.
If we can express
this difference as a linear combination of elements $p$ 
in the ideal $(P)$ with $\LM(p)<m$, 
then we say the ambiguity was resolvable {\em relative to the ordering}.
(Bokut' used ``closed under composition'' to describe this
condition along with minimality of $P$.)
See~\cite[Theorem~1.2]{Bergman}.
%%%%%%%%%%%%%%%%%%%%%%%%%%%%%%%%%%%%%%
\begin{namedthm}[Diamond Lemma idea]
The following are equivalent:
\begin{itemize}
\item
The set of reduced words is a $k$-basis of $T/(P)$.
\item All ambiguities among elements of $P$ are resolvable.
\item All ambiguities among elements of $P$ are resolvable relative to
the ordering.
\item Every element in $(P)$ can be reduced to zero in $T/(P)$
by just using the relations in $P$.
\end{itemize}
\end{namedthm}
%%%%%%%%%%%%%%%%%%%%%%%%%%%%%%%%%%%%%%%%%%%%%%%%%
In essence, the lemma says that if the generating set of relations
$P$ is well-behaved with respect to some monomial ordering,
then one can define a canonical form just by checking that
nothing goes wrong with the set $P$ instead of checking
for problems implied by the whole ideal $(P)$.
Thus, resolving ambiguities is just
another way of testing for a Gr\"obner basis
(see~\cite{Green94}):
%%%%%%%%%%%%%%%%%%%%%%%%%%%%%%%%%%%%%%%%%%%%%%%%%%%%%%
The set $P$ is a Gr\"obner basis for the ideal $(P)$ 
if and only if all ambiguities among elements of $P$ are resolvable.
%%%%%%%%%%%%%%%%%%%%%%%%%%%%%%%%%%%%%%%%%%%%%%%%%%

%%%%%%%%%%%%%%%%%%%%%%%%%%%%%%%%%%%%%%%%%%%%%%%5
\subsection*{Applications}
Although the idea of the Composition-Diamond lemma
can be phrased in many ways, the hypothesis to be checked
in the various versions of the lemma requires very similar computations
in application. One finds precursors
of the ideas underlying the Composition-Diamond Lemma in the original proofs
given by Poincar\'e, Birkhoff, and Witt
of the PBW Theorem for universal
enveloping algebras of Lie algebras.
These techniques and computations have been used in a number of other settings.
For example, explicit
PBW conditions are given
for Drinfeld Hecke algebras (which include symplectic reflection
algebras) by Ram and Shepler~\cite{RamShepler}; see Section~\ref{SRA}.
In~\cite{doa}, 
we studied the general  case of algebras defined 
by relations which set commutators to lower order terms
using both a homological approach
and the Composition-Diamond Lemma 
(as it holds 
in arbitrary characteristic).
These algebras, called {\em Drinfeld orbifold algebras},
include 
Sridharan enveloping algebras, Drinfeld Hecke algebras, enveloping algebras of Lie algebras,
Weyl algebras, and twists of these algebras with group actions.
Gr\"obner bases were used explicitly by Levandovskky 
and Shepler~\cite{LevandovskyyShepler}
in replacing a commutative polynomial algebra by a skew (or quantum)
polynomial algebra in the theory of Drinfeld Hecke algebras.
Bergman's Diamond Lemma was used by Khare~\cite{Khare} to generalize
the Drinfeld Hecke algebras of Section~\ref{SRA} from the setting of group actions
to that of algebra actions.

Of course the Composition-Diamond Lemma and Gr\"obner-Shirshov bases
have been used to explore many different kinds of algebras
(and in particular to find PBW-like bases) that we will 
not discuss here. 
See Bokut' and Kukin~\cite{BokutKukin} and
Bokut' and Chen~\cite{BokutChen} for many such examples.

Note that some authors prove PBW theorems by creating a space
upon which the algebra in question acts
(see, e.g., Humphreys~\cite{Humphreys} or Griffeth~\cite[first version]{Griffeth}).  
Showing that the given space is
actually a module for the algebra requires checking certain relations
that are similar to the conditions that one must check before invoking
the Composition-Diamond Lemma.

%%%%%%%%%%%%%%%%%%%%%%%%%%%%%%%%%%%%%%%%%%%%%%%%%%%%%%%%%%%%
%%%%%%%%%%%%%%%%%%%%%%%%%%%%%%%%%%%%%%%%%%%%%%%%%%%%%%%%%%%%
%%%%%%%%%%%%%%%%%%%%%%%%%%%%%%%%%%%%%%%%%%%%%%%%%%%%%%%%%%%%
\section{Drinfeld-Jimbo quantum groups and related Hopf algebras}\label{DJQG} 

Quantized  enveloping algebras (that is, Drinfeld-Jimbo quantum
groups~\cite{Drinfeld87,Jimbo})
are deformations of universal enveloping algebras of Lie algebras.
(Technically, they are bialgebra deformations rather than algebra deformations.)
Many mathematicians discovered PBW bases for these algebras,
in particular Lusztig~\cite{Lusztig1,Lusztig1.5,Lusztig2} 
in a very general setting
and DeConcini and Kac~\cite{DeConciniKac} 
by defining a corresponding algebra filtration.
Although 
there are many incarnations of these algebras, 
we restrict ourselves to 
the simply-laced case
and to algebras over the complex numbers for ease of notation. 
We state a PBW theorem in this context and 
refer the reader to the literature for more general 
statements (see, e.g.,~\cite{Lusztig2}).   

\subsection*{Quantum groups}
Let $\mathfrak g$ be 
a finite dimensional semisimple complex Lie algebra  of rank $n$
with symmetric Cartan matrix $(a_{ij})$.
Let $q$ be a nonzero complex number, $q\neq \pm 1$.
(Often $q$ is taken to be an indeterminate instead.) 
The {\em quantum group} $U_q({\mathfrak g})$ is the associative $\CC$-algebra
defined by generators
$$E_1,\ldots, E_n, F_1,\ldots,F_n, K_1^{\pm 1} ,\ldots, K_n^{\pm 1}$$ and
relations
\begin{eqnarray*}
  K_i ^{\pm 1} K_j^{\pm 1}   =  K_j^{\pm 1}K_i^{\pm 1}, &  &
    K_i K_i^{-1}  =  1   =   K_i^{-1}K_i,\\
  K_i E_j  =   q^{ a_{ij}} E_j K_i , &  & K_i F_j
   =  q^{- a_{ij}} F_j K_i, \\
  E_i F_j - F_j E_i \! & = & \! \delta_{ij} \, \frac{K_i - K_i^{-1}}{q-q^{-1}},\\
 E_i^2 E_j - (q+q^{-1}) E_i E_j E_i + E_j E_i^2\!  & = & \!  0  \ \ \mbox{ if }a_{ij}=-1,
    \ \ \ E_iE_j=E_jE_i \ \ \mbox{ if }a_{ij}=0,\\
 F_i^2 F_j \, - (q+q^{-1})\, F_iF_jF_i \, + F_j F_i^2 \! & = &\!  0 \ \ \mbox{ if }a_{ij}=-1,
   \ \ \ F_iF_j\, = F_jF_i \  \ \ \mbox{ if }a_{ij}=0.  
\end{eqnarray*}
The last two sets of relations are called the quantum Serre relations. 

Let $W$ be the Weyl group of $\mathfrak g$. Fix a reduced
expression of the longest element $w_0$ of $W$. 
This choice yields a total order on the set $\Phi^+$ of positive roots,
$ \ \beta_1,\ldots, \beta_m$. 
To each $\alpha\in\Phi^+$, 
Lusztig~\cite{Lusztig1,Lusztig1.5,Lusztig2} assigned an  element 
$E_{\alpha}$ (respectively, $F_{\alpha}$) in $U_q({\mathfrak {g}})$
determined by this ordering 
that is an iterated braided commutator of the generators
$E_1,\ldots, E_n$ (respectively, $F_1,\ldots, F_n$).
These ``root vectors'' then appear in a PBW basis: 

%%%%%%%%%%%%%%%%%%%%%%%%%%%%%%%%%%%%%%%%%%%
\begin{namedthm}[PBW Theorem for Quantum Groups]
There  is a basis of $U_q({\mathfrak{g}})$ given by 
$$
  \{ E_{\beta_1}^{a_1}\cdots E_{\beta_m}^{a_m}
     K_1^{b_1}\cdots K_n^{b_n} 
    F_{\beta_1}^{c_1}\cdots F_{\beta_m}^{c_m}  : 
     a_i, c_i \geq 0, \ b_i\in \Z \}.
$$
Moreover, there is a filtration on the subalgebra $ U_q^{>0}({\mathfrak{g}})$ (respectively,
$U_q^{<0}({\mathfrak{g}})$) generated by $E_1,\ldots,E_n$ (respectively,
$F_1,\ldots,F_n$) for which the  associated graded algebra 
is isomorphic to a skew polynomial ring.
\end{namedthm}
%%%%%%%%%%%%%%%%%%%%%%%%%%%%%%%%%%%%%%%%%%%%

The skew polynomial ring to which the theorem refers is generated by the
images of the $E_{\alpha}$ (respectively, $F_{\alpha}$), with relations
$E_{\alpha}E_{\beta} = q_{\alpha\beta}E_{\beta}E_{\alpha}$
(respectively, $F_{\alpha}F_{\beta} = q_{\alpha\beta}F_{\beta}F_{\alpha}$) 
where each 
$q_{\alpha\beta}$
is a scalar determined by $q$ and by $\alpha,\beta$ in $\Phi^+$. 

\vspace{2ex}
\begin{ex}
The algebra $U^{>0}_q({\mathfrak{sl}}_3)$ is generated by $E_1,E_2$.
Let $$E_{12}:= E_1E_2 - qE_2E_1.$$
Then, as a consequence of the quantum Serre relations, $E_1E_{12}= q^{-1}E_{12}E_1$
and 
$E_{12}E_2=q^{-1}E_2E_{12}$, and, by definition of $E_{12}$, we also have 
$E_1E_2= qE_2E_1 + E_{12}$.
In the associated graded algebra, this last relation becomes $E_1E_2=qE_2E_1$. 
The algebras $U^{>0}_q({\mathfrak{sl}}_n)$ are similar, however in
general the filtration on $U^{>0}_q({\mathfrak{g}})$ stated in
the theorem is more complicated.
\end{ex}
\vspace{2ex}

\subsection*{Proofs and related results}
There are several proofs in the literature of the first statement of 
the above theorem
 and related results,
beginning with Khoroshkin and Tolstoy~\cite{KhoroshkinTolstoy},
Lusztig~\cite{Lusztig1,Lusztig1.5,Lusztig2}, Takeuchi~\cite{Takeuchi},
and Yamane~\cite{Yamane}.
These generally involve explicit computations facilitated by 
representation theory.
Specifically, one obtains representations of $U_q({\mathfrak{g}})$ from those of the
corresponding Lie algebra $\mathfrak g$ by deformation, and one
then uses what is known in the classical setting to obtain
information about $U_q({\mathfrak{g}})$. 
Ringel~\cite{Ringel} gave a different approach via Hall algebras.
The filtration  and structure of the associated graded algebra
of $U^{>0}(\mathfrak{g})$  
was first given by DeConcini and Kac~\cite{DeConciniKac}.
For a general ``quantum PBW Theorem'' that applies to some of these
algebras, see Berger~\cite{Berger}. 

In case $q$ is a root of unity (of order $\ell$), there are
finite dimensional versions of Drinfeld-Jimbo quantum groups.
The {\em small quantum group} $u_q({\mathfrak{g}})$ may be defined
as the quotient of $U_q({\mathfrak{g}})$ by the ideal generated by all
$E^{\ell}_{\alpha}, F^{\ell}_{\alpha}, K^{\ell}_{\alpha} - 1$.
This finite dimensional algebra has $k$-basis given by elements in
the PBW basis of the above theorem for which
$0\leq a_i,b_i,c_i <\ell$.

The existence of 
PBW bases for $U_q({\mathfrak{g}})$ and $u_q({\mathfrak{g}})$ 
plays a crucial role in their representation theory,
just as it does in the classical setting of Lie algebras. 
Bases of finite dimensional simple modules and other modules are defined from
weight vectors and PBW bases~\cite{Lusztig1.5}.
R-matrices may be expressed 
in terms of PBW basis elements~\cite{Drinfeld87,Jimbo,Rosso}.
Computations of cohomology take advantage of the structure provided by
the PBW basis and filtration (e.g., see~\cite{GinzburgKumar}, based on techniques developed
for restricted Lie algebras~\cite{FriedlanderParshall}). 

More generally, PBW bases and some Lie-theoretic structure appear in a
much larger class of Hopf algebras. Efforts to understand 
finite dimensional
Hopf algebras of various types led in particular to a study of those arising from
underlying Nichols algebras. 
Consequently, a classification of some types of pointed Hopf algebras was completed by 
Andruskiewitsch
and Schneider~\cite{AndruskiewitschSchneider},
Heckenberger~\cite{Heckenberger}, and Rosso~\cite{Rosso2}. 
A Nichols algebra is a  ``braided'' graded Hopf algebra that is 
connected, generated by its 
degree~1 elements, and whose  subspace of primitive elements
is precisely its degree~1 component. 
The simplest Nichols algebras are those of ``diagonal type,''
and these underlie the Drinfeld-Jimbo quantum groups and the 
Hopf algebras in the above-mentioned  classification. These
algebras have
PBW bases just as does $U^{>0}_q({\mathfrak{g}})$ or $u_q^{>0}
({\mathfrak{g}})$;
a proof given by Kharchenko~\cite{Kharchenko} uses a combinatorial 
approach such as in Section~\ref{diamond}.

%%%%%%%%%%%%%%%%%%%%%%%%%%%%%%%%%%%%%%%%%%%%%%%%%%%%%%%%%%%%
%%%%%%%%%%%%%%%%%%%%%%%%%%%%%%%%%%%%%%%%%%%%%%%%%%%%%%%%%%%%
%%%%%%%%%%%%%%%%%%%%%%%%%%%%%%%%%%%%%%%%%%%%%%%%%%%%%%%%%%%%
\section{Symplectic reflection algebras, rational
Cherednik algebras, and graded (Drinfeld) Hecke algebras}\label{SRA}

Drinfeld~\cite{Drinfeld86} and Lusztig~\cite{Lusztig88, Lusztig89} originally defined
the algebras now variously called symplectic reflection algebras, rational
Cherednik algebras, and graded (Drinfeld) Hecke algebras, depending on context.
These are PBW deformations of group extensions of polynomial rings
(skew group algebras) defined by relations
that set commutators of vectors to elements of a group algebra. 
Lusztig explored the representation theory of
these algebras when the acting group is a Weyl group.  
Crawley-Boevey and Holland~\cite{CBH} considered subgroups of ${\rm SL}_2(\CC)$ 
and studied subalgebras of these algebras in relation to corresponding orbifolds. 
Initial work on these types of PBW deformations for arbitrary groups 
began with
Etingof and Ginzburg~\cite{EtingofGinzburg} and Ram and Shepler~\cite{RamShepler}. 
Gordon~\cite{GordonInvent} used the rational Cherednik algebra to
prove a version of the $n !$-conjecture for Weyl groups 
and the representation theory of these algebras remains
an active area.  
(See~\cite{BrownSurvey},~\cite{GordonSurveyCherednik},~\cite{GordonSurveySymplectic},
and~\cite{RouquierSurvey}.)  
We briefly recall and compare these algebras.
(See also
~\cite{Chlouveraki} 
for a survey of symplectic reflection algebras
and rational Cherednik algebras in the context of Hecke algebras
and representation theory.) 

Let $G$ be a group acting by automorphisms on a $k$-algebra $S$. 
The {\em skew group algebra} $S\# G$ (also written as a semidirect y
product $S\rtimes G$)
 is the $k$-vector space $S\otimes kG$
together with multiplication given by 
$(r \otimes g)(s\otimes h)=r ( \,^{g}s)\otimes g h$
for all $r,s$ in $S$ and $g,h$ in $G$,
where $^gs$ is the image of $s$ under the automorphism $g$.

%%%%%%%%%%%%%%%%%%%%%%%%%%%%%%%%%%%%
\subsection*{Drinfeld's ``Hecke algebra''}
Suppose $G$ is a finite group acting linearly on a finite dimensional
vector space $V$ over $k=\CC$ with symmetric algebra $S(V)$.
Consider the quotient algebra
$$
\cH_{\kappa}=\quotient{T(V)\#G}{( v_1\otimes v_2-v_2\otimes v_1-\kappa(v_1,v_2): 
v_1,v_2\in V) } $$
defined by a bilinear
parameter function $\kappa:V\times V\rightarrow \CC G$.
We view $\cH_{\kappa}$ as a filtered algebra by assigning degree one to
vectors in $V$ and degree zero to group elements in $G$.
Then the algebra $\cH_{\kappa}$ is a PBW deformation of $S(V)\# G$ 
if its associated graded algebra is isomorphic to $S(V)\# G$.
Drinfeld~\cite{Drinfeld86} originally
defined these algebras for arbitrary groups, and 
he also 
counted the dimension of the parameter space
of such PBW deformations for Coxeter groups (see also~\cite{RamShepler}).

\vspace{2ex}
\begin{ex} 
Let $V$ be a vector space of dimension 3 with basis $v_1,v_2,v_3$,
and let $G$ be the symmetric group $S_3$ acting on $V$
by permuting the chosen basis elements. 
The following is a PBW deformation of $S(V)\# G$, where
$(i\, j\, k)$ denotes a 3-cycle in $S_3$: 
$$\cH_{\kappa} = \quotient{T(V)\# S_3}
   { ( v_i\otimes v_j - v_j\otimes v_i - (i\, j\, k) + (j\, i\, k) :
    \{i,j,k\}=\{1,2,3\} ) } . 
$$
\end{ex}
\vspace{2ex}

%%%%%%%%%%%%%%%%%%%%%%%%%%%%%%%5
\subsection*{Lusztig's graded affine Hecke algebra}
While exploring the representation theory of groups of Lie type,
Lusztig~\cite{Lusztig88, Lusztig89}
defined a variant 
of the affine Hecke algebra for Weyl groups
which he called ``graded''
(as it was obtained from a particular filtration
of the affine Hecke algebra).
He gave a presentation for this algebra $\mathbb{H}_{\lambda}$ using the
same generators as those for Drinfeld's Hecke algebra
$\cH_{\kappa}$, but he gave relations  preserving 
the structure of the polynomial ring 
and altering the skew group algebra relation.
(Drinfeld's relations do the reverse.)
The {\em graded affine Hecke algebra} $\mathbb{H}_{\lambda}$ 
(or simply the ``graded Hecke algebra'')
for a finite Coxeter group $G$ acting on a finite dimensional complex
vector space $V$ (in its natural reflection representation) is
the $\CC$-algebra generated by the polynomial algebra $S(V)$ 
together with the group algebra
$\CC G$ with relations $$g v =\ ^g v g+ \lambda_g(v) g$$
for all $v$ in $V$ and $g$ in a set $\mathcal S$ of simple reflections
(generating $G$) where $\lambda_g$ in $V^*$ defines the reflecting hyperplane
($\text{ker}\, \lambda_g\subseteq V$) 
of $g$ and $\lambda_g=\lambda_{hgh^{-1}}$ for all $h$ in $G$.
(Recall that a {\em reflection}
on a finite dimensional vector space is just
a nonidentity linear transformation 
that fixes a hyperplane pointwise.) 

Note that for $g$ representing a fixed conjugacy class of reflections,
the linear form $\lambda_g$ is only well-defined up to a
nonzero scalar.  Thus one often fixes once and for all
a choice of linear forms $\lambda = \{ \lambda_g \}$
defining the orbits of reflecting hyperplanes
(usually expressed using Demazure/BGG operators) and then
introduces a formal parameter by which to rescale.
This highlights the degree of freedom arising from each orbit; 
for example, one might replace
$$
\lambda_g(v)\ \ \text{ by }\ \ 
c_g \ \langle v, \alpha_g^{\vee} \rangle 
= c_g\ \Big( \frac{v-\, ^gv}{\alpha_g}\Big)
$$
for some conjugation invariant formal
parameter $c_g$ after fixing a 
$G$-invariant inner product and root system $\{\alpha_g:
g\in \mathcal S\} \subset V$
with coroot vectors $\alpha_g^{\vee}$.
(Note that for any reflection $g$,
the vector $(v-\, ^gv)$ is a nonzero scalar multiple of $\alpha_g$
and so the quotient of $v-\, ^gv$ by $\alpha_g$ is a scalar.)
Each graded affine Hecke algebra $\mathbb{H}_{\lambda}$
is filtered 
with vectors in degree one and group elements in degree zero
and defines a PBW deformation of $S(V)\# G$.

\vspace{2ex}

\begin{center}
\begin{fbox}
{\begin{tabular}{|l | l | l |}
\hline
$\rule[-1ex]{0ex}{4ex}$
Finite Group  & Any $G\leq \text{GL}(V)$ & 
Coxeter $G \leq \text{GL}(V)$\\
\hline
$\rule[-1ex]{0ex}{4ex}$
Algebra
& $\cH_{\kappa}$  \text{(Drinfeld)}
& $\mathbb{H}_{\lambda}$ \text{(Lusztig)}\\
\hline
$\rule[-1ex]{0ex}{4ex}$
generated by & $V\text{ and } \CC G$ & $ V\text{ and } \CC G$ \\
\hline
$\rule[-1ex]{0ex}{4ex}$
with relations & $ g v =\ ^g v g , $ &  $g v =\, ^g vg + \lambda_g(v)g ,$\\
$\rule[-1ex]{0ex}{2ex}$
& $v w = w v + \kappa(v,w) $  & $v w = w v $\\
& $(\forall v,w\in V, \ \forall g\in G)$ & $(\forall v,w\in V, \ \forall g\in \mathcal S)$\\
\hline
\end{tabular}}
\end{fbox}
\end{center}

\vspace{2ex}

\subsection*{Comparing algebras}
Ram and Shepler~\cite{RamShepler} showed that Lusztig's graded
affine Hecke algebras are a special case of Drinfeld's construction:
For each parameter $\lambda$, there is a parameter $\kappa$
so that the filtered algebras
$\mathbb{H}_{\lambda}$ and $\cH_{\kappa}$ are isomorphic (see~\cite{RamShepler}).
Etingof and Ginzburg~\cite{EtingofGinzburg} rediscovered 
Drinfeld's algebras with focus on groups $G$ acting symplectically 
(in the context of orbifold theory).  They called algebras $\cH_{\kappa}$
satisfying the PBW property {\em symplectic reflection algebras},
giving necessary
and sufficient conditions on $\kappa$ for symplectic groups.
They used the theory of 
Beilinson, Ginzburg, and Soergel~\cite{BGS} of Koszul rings
to generalize Braverman and Gaitsgory's conditions
to the setting where the ground field  is replaced by the semisimple 
group ring $\CC G$.  (The skew group algebra $S(V)\# G$ is Koszul
as a ring over the semisimple subring $\CC G$.)
Ram and Shepler~\cite{RamShepler}
independently gave necessary and sufficient PBW 
conditions on $\kappa$ for arbitrary groups acting
linearly over $\CC$
and classified all such quotient 
algebras for complex reflection groups. Their proof 
relies on the Composition-Diamond Lemma.
(See Sections~\ref{BG} and~\ref{diamond} for
a comparison of these two techniques for showing PBW properties.)
Both approaches depend on the fact that the underlying field $k=\CC$
has characteristic zero (or, more generally, has characteristic
that does not divide the order of the group $G$).  
See Section~\ref{positivechar} for a discussion of PBW theorems
in the modular setting when $\CC$ is replaced
by a field whose characteristic divides $|G|$.

%%%%%%%%%%%%%%%%%%%%%%%%%%%%%
\subsection*{Rational Cherednik algebras}
The rational Cherednik algebra is a special case of a quotient
algebra $\cH_{\kappa}$ satisfying the PBW property
(in fact,
a special case of a symplectic reflection algebra) for 
reflection groups acting diagonally on two copies
of their reflection representations (or ``doubled up'').
These algebras are regarded as ``doubly degenerate'' versions of the
double affine Hecke algebra
introduced by Cherednik~\cite{Cherednik}
to solve the Macdonald
(constant term) conjectures in combinatorics.  We simply
recall the definition here in terms of reflections
and hyperplane arrangements.

Suppose $G$ is a finite group generated by reflections on 
a finite dimensional complex vector space $V$.  
(If $G$ is a Coxeter group, then extend the action to
one over the complex numbers.)
Then the induced diagonal action
of $G$ on $V\oplus V^*$ is generated by {\em bireflections}
(linear transformations that each fix a subspace
of codimension 2 pointwise), i.e., by {\em symplectic reflections} 
with respect to a natural symplectic form on $V\oplus V^*$.

Let $\mathcal{R}$ be the set of all
reflections in $G$ acting on $V$. 
For each reflection $s$ in $\mathcal{R}$,
let
$\alpha_s$ in $V$ and $\alpha_s^*$ in $V^*$ 
be eigenvectors (``root vectors'') each
with nonindentity eigenvalue.
We define 
an algebra generated by $\CC G$, $V$, and $V^*$
in which vectors in $V$ commute with each other and vectors in $V^*$
commute with each other, 
but passing a vector from $V$ over one from $V^*$ gives
a linear combination of reflections (and the identity).
As parameters, we take a scalar
$t$ and a set of scalars ${\bf c}=\{c_s:s\in \mathcal{R}\}$ 
with $c_s=c_{hsh^{-1}}$ for all $h$ in $G$.
The {\em rational Cherednik algebra} ${\text{\bf H}}_{t,{\bf c}}$
with parameters
$t, {\bf c}$ 
is then the $\CC$-algebra generated by the vectors in $V$ and $V^*$
together with the group algebra $\CC G$ 
satisfying the relations 
$$
\begin{aligned}
gu & =\ ^g u g, \quad
 u u' =u' u,\\ 
v v^*& = v^* v+t\, v^*(v)
-\sum_{s\in \mathcal{R}} c_s\ \alpha_s^*(v)\ v^*(\alpha_s)\ s
\end{aligned}
$$
for any $g$ in $G$, $v$ in $V$, $v^*$ in $V^*$,
and any $u,u'$ both in $V$ or both in $V^*$.
Note that $\alpha_s$ and $\alpha_s^*$
are only well-defined up to a nonzero scalar,
and we make some conjugation invariant choice of
normalization in this definition, say, 
by assuming that $\alpha_s^*(\alpha_s)=1$.
One often replaces $\CC$ by  $\CC[t, {\bf c}]$ to work
in a formal parameter space. 

The relations defining the rational Cherednik algebra are often 
given in terms of the 
arrangement of reflecting hyperplanes $\mathcal{A}$ for $G$ acting on $V$.
For each hyperplane $H$ in $\mathcal{A}$, 
choose a linear form $\alpha_H^*$ in $V^*$ defining
$H$ (so $H=\text{ker}\, \alpha_H^*$) and let $\alpha_H$  be a nonzero vector
in $V$ perpendicular to $H$ with respect to some fixed
$G$-invariant inner product.  Then the third defining relation of
${\text{\bf H}}_{t,{\bf c}}$ 
can be rewritten (without a choice of normalization) as 
$$v v^*-v^* v= t v^*(v) -
\sum_{H\in \mathcal{A}}
\frac{\alpha_H^*(v)\ v^*(\alpha_H)}
{\alpha_H^*(\alpha_H)}
\big(c_{s_H} s_H + c_{s_H^2} s_H^2 + \ldots +c_{s_H^{a_H}} s_H^{a_H}\big)
$$
where $s_H$ is the reflection in $G$ about 
the hyperplane $H$ of maximal order
$a_H+1$.
Again, this is usually expressed geometrically in terms of 
the inner product on $V$ and induced product on $V^*$:
$$\frac{\alpha_H^*(v)\ v^*(\alpha_H)}
{\alpha_H^*(\alpha_H)} = 
\frac{\langle v, \alpha_H^{\vee}\rangle
\langle \alpha_H, v^*\rangle}
{\langle \alpha_H, \alpha^{\vee}\rangle }\ .
$$

The PBW theorem then holds for the algebra ${\text{\bf H}}_{t,{\bf c}}$
(see~\cite{EtingofGinzburg}):
%%%%%%%%%%%%%%%%%%%%%%%%%%%%%%%%%%%%%%%%%%%%%%%%%%%
\begin{namedthm}[PBW Theorem for Rational Cherednik Algebras]
The rational Cherednik algebra
${\text{\bf H}}_{t,{\bf c}}$ is isomorphic to $S(V)\otimes S(V^*)\otimes \CC G$
as a complex vector space for any choices of parameters $t$ and $c$,
and its associated graded
algebra is isomorphic to $(S(V)\otimes S(V^*))\# G$.
\end{namedthm}
%%%%%%%%%%%%%%%%%%%%%%%%%%%%%%%%%%%%%%%%%%%%%%%%%%%%%%%%%%%%%

Connections between rational Cherednik algebras 
and other fields of mathematics are growing stronger.
For example, Gordon and Griffeth~\cite{GordonGriffeth} link the
Fuss-Catalan numbers in combinatorics to the representation theory
of rational Cherednik algebras.
These investigations also bring insight
to the classical theory of complex
reflection groups,
especially to the perplexing
question of why some reflection groups
acting on $n$-dimensional space
can be generated by $n$ reflections (called ``well-generated'' or 
``duality'' groups) and others not.
(See~\cite{BerkeschGriffethSam, 
GorskyOblomkovRasmussenShende,
ShanVaragnoloVasserot}
for other recent applications.)

%%%%%%%%%%%%%%%%%%%%%%%%%%%%%%%%%%%%%%%%%%%%%%%%%%%%%%%%%%%%%%%%%%%
%%%%%%%%%%%%%%%%%%%%%%%%%%%%%%%%%%%%%%%%%%%%%%%%%%%%%%%%%%%%%%%%%%%
%%%%%%%%%%%%%%%%%%%%%%%%%%%%%%%%%%%%%%%%%%%%%%%%%%%%%%%%%%%%%%%%%%%
\section{Positive characteristic and nonsemisimple ground rings}\label{positivechar}
Algebras displaying PBW properties are quite common over ground
fields of positive characteristic and nonsemisimple ground rings, but
techniques for establishing PBW theorems are not all
equally suited for work over arbitrary fields and rings.
We briefly mention a few results of ongoing efforts to establish 
and apply PBW theorems
in these  settings.

The algebras of Section~\ref{SRA} make sense in the {\em modular setting},
that is, when the characteristic of $k$ is a prime dividing the order of
the finite group $G$. In this case, however, the group algebra $kG$ is
not semisimple, and one must take more care in proofs.  
PBW conditions on $\kappa$ were examined by
Griffeth~\cite{Griffeth}
by construction of an explicit
$\cH_{\kappa}$-module, as is done in one standard proof of the PBW 
Theorem for universal enveloping algebras. 
(See also Bazlov and
Berenstein~\cite{BazlovBerenstein} for a generalization.) 
The Composition-Diamond Lemma, being characteristic free, applies in the modular setting;
see our paper~\cite{doa} for a proof of the PBW property using
this lemma that applies to graded
(Drinfeld) Hecke algebras over fields of arbitrary characteristic.
(Gr\"obner bases are explicitly used 
in Levandovskyy and Shepler~\cite{LevandovskyyShepler}.)
Several authors consider
representations of rational Cherednik algebras in the modular setting,
for example, Balagovic and Chen~\cite{BalagovicChen},
Griffeth~\cite{Griffeth}, 
and Norton~\cite{Norton}.

The theory of Beilinson, Ginzburg, and Soergel of Koszul rings over
semisimple subrings, used in Braverman-Gaitsgory style proofs of PBW theorems,
does not apply directly to the modular setting. 
However it may be adapted using 
a larger complex replacing  the Koszul complex:  
In~\cite{PBWQuadratic}, we used this approach to
 generalize the Braverman-Gaitsgory argument
to arbitrary Koszul algebras with finite group actions. 
This replacement complex has an advantage over the Composition-Diamond Lemma
or Gr\"obner basis theory arguments in that it contains information about 
potentially new types of deformations that do not occur in the nonmodular setting.

Other constructions generalize the
algebras of Section~\ref{SRA} to algebras over ground rings 
that are not necessarily
semisimple. Etingof, Gan, and Ginzburg~\cite{EGG} considered 
deformations of algebras that are extensions of polynomial rings by 
acting algebraic groups or Lie algebras.  They
used a Braverman-Gaitsgory approach to obtain a Jacobi condition
by realizing the acting algebras as inverse limits of finite dimensional
semsimple algebras. 
Gan and Khare~\cite{GanKhare}
investigated actions of $U_q({\mathfrak{sl}}_2)$ on the quantum
plane (a skew polynomial algebra), and Khare~\cite{Khare}
looked at actions of arbitrary cocommutative
algebras on polynomial rings.
In both cases  PBW theorems were proven using the Composition-Diamond 
Lemma. 
A general result for actions of (not necessarily semisimple) Hopf
algebras on Koszul algebras is contained in 
Walton and Witherspoon~\cite{WaltonWitherspoon}
with a Braverman-Gaitsgory style proof. 
See also He, Van Oystaeyen, and Zhang~\cite{HVZ} for a PBW theorem using a
somewhat different complex in a general setting of Koszul rings over  
not necessarily semisimple ground rings. 
One expects yet further generalizations and applications of 
the ubiquitous and 
potent 
PBW Theorem.

%%%%%%%%%%%%%%%%%%%%%%%%%%%%%%%%%%%%%%%%%%%%%%%%%%%%%%%%%%%%%
%%%%%%%%%%%%%%%%%%%%%%%%%%%%%%%%%%%%%%%%%%%%%%%%%%%%%%%%%%%%%
%%%%%%%%%%%%%%%%%%%%%%%%%%%%%%%%%%%%%%%%%%%%%%%%%%%%%%%%%%%%

\end{document}